\documentclass[12pt,oneside,english]{amsart}
\usepackage[latin1]{inputenc}
\usepackage{amssymb}

\makeatletter

\newtheorem{lemma}{Lemma}

\usepackage{babel}

\makeatother
\begin{document}
\baselineskip=17pt

\title{On excess of the odious primes}

\author{Vladimir Shevelev}
\address{Departments of Mathematics \\Ben-Gurion University of the
 Negev\\Beer-Sheva 84105, Israel. e-mail:shevelev@bgu.ac.il}

\subjclass{11N05.}

\begin{abstract}
We give a more strong heuristic justification of our conjecture on
the excess of the odious primes.
\end{abstract}

\maketitle
\section{Introduction}
     This note is a continuation of the author's paper \cite{4}.
Until now the author did not know about the Moser's digit conjecture
and its solutions in \cite{3} and \cite{1}. In fact, we gave in
\cite{4} a new combinatorial proof of this conjecture  (see Theorem
3\cite{4}) and proved also an addition to the Moser-Newman theorem
for the excess of the evil nonnegative odd integers less than $n$
and divisible by 3(see. Theorem 4 and 5 in \cite{4}).

     The aim of the present note - to give a more strong heuristic
justification ("almost strong proof") of our Conjectures 1 and 2
\cite{4}.

     Recently in their excellent research  \cite{2}, Mauduit and Rivat  solved the
Gelfond digit problem for primes. In particular, they proved that

\begin{equation}\label{1}
\lim_{n\rightarrow\infty}\frac{\pi^0(n)}{\pi(n)}=\lim_{n\rightarrow\infty}
\frac{\pi^e(n)}{\pi(n)}=\frac 1 2.
\end{equation}

Moreover, if $\pi^0_{3,i}(n)(\pi^e_{3,i}(n)),\;\;i=1,2$, is the
number of the odious (evil) primes $p\equiv i(mod{3})$ not exceeding
$n$ and

$$
\pi_{3,i}(n)=\pi^0_{3,i}(n)+\pi^e_{3,i}(n),\;\;i=1,2,
$$

they also proved that

\begin{equation}\label{7}
\lim_{n\rightarrow\infty}\frac{\pi^0_{3,i}(n)}{\pi_{3,i}(n)}=
\lim_{n\rightarrow\infty}\frac{\pi^e_{3,i}(n)}{\pi_{3,i}(n)}=\frac 1
2
\end{equation}

 These results mean that the events "$n$ is a prime" and "$n$ is an
odious integer" are asymptotic independent for large $n$.

In turn, this means that the odious-evil asymptotic behavior of the
primes of the form $3k+1(3k+2)$ is proportionally similar to the
odious-evil asymptotic behavior of all odd integers of the form
$3k+1(3k+2)$.

\newpage

\section{Proof of conjectures}

     Let $\mu^0(n)(\mu^e(n))$ be the number of \slshape odd \upshape odious (evil)
integers less than $n$.

\begin{lemma}      $|\mu^0(n)-\mu^e(n)|\leq 1$.
\end {lemma}

\slshape Proof. \upshape The lemma follows from the identity

\begin{equation}\label{6}
\mu^0(4m+1)=\mu^e(4m+1),\;\;m\in\mathbb{N},
\end{equation}

which is proved by induction.

     Notice that (\ref{6}) is valid for $m=1$. Assuming that it is
valid for $4m+1$ we prove (\ref{6}) for $4(m+1)+1$. Indeed, let $m$
have $k$ ones in the binary expansion. Then taking into account that
for odd $k$ the number $4m+1$ is evil and for even $k$ the number
$4m+1$ is odious, and using the induction conjecture we have

$$
\mu^0(4m+3)-\mu^e(4m+3)=(-1)^k.
$$

     Furthermore, $4m+3$ is evil if $k$ is even and is odious if $k$
is odd. Therefore,

$$
\mu^0(4m+5)-\mu^e(4m+5)=(-1)^k+(-1)^{k+1}=0. \blacksquare
$$

     Let $\Delta_{3,i}(n)(\Delta^{odd}_{3,i}(n)),\;\;i=0,1,2$  be the excess of the
\slshape (odd )\upshape  odious integers $m\in[0,n)$ such that
$m\equiv i \;(\mod{3})$.

In particular, according to the notations of \cite{4}

\begin{equation}\label{2}
\Delta^{odd}_{3,0}(n)=-\Delta^{odd}_3([0,n))< 0,\;\;\Delta
_{3,0}(n)=-\Delta_3([0,n))<0.
\end{equation}

     Let, furthermore, $\Delta^{primes}_{3,i},\;\;i=1,2$  be the
excess of odious odd primes $p\in[0,n)$ such that $p\equiv
i(\mod{3})$.  Then by (\ref{1}),( \ref{7}) taking into account the
independence of the above mentioned events, in the case of
$|\Delta^{odd}_{3,i}(n)|>>\ln{n}$ we have

\begin{equation}\label{3}
\Delta^{primes}_{3,i}\sim\frac{3\Delta^{odd}_{3,i}(n)}{\ln
n},\;\;\;i=1,2.
\end{equation}

So, for $i=1$ and for even $n$ we have

\begin{equation}\label{8}
\Delta^{odd}_{3,1}(n)=-\Delta_{3,0}\left(\frac n
2\right)=\Delta_3\left([0.\frac n 2)\right).
\end{equation}
\newpage

Newman showed \cite{3}, that for all $n\in\mathbb{N}$

\begin{equation}\label{9}
(0.05\cdot 3^\alpha)n^\alpha\leq \Delta_3([0,n))\leq(5\cdot
3^\alpha)n^\alpha\;\;\;with \;\;\alpha=\frac{\ln{3}}{\ln{4}}.
\end{equation}

Therefore, by (\ref{8})

\begin{equation}\label{10}
\Delta^{odd}_{3,1}(n)\geq 0.05(1.5)^\alpha n^\alpha >> \ln{n}.
\end{equation}

In the case of $i=2$ the absolute value of the excess
$\Delta^{odd}_{3,2}(n)$  is small for some $n$. Indeed, by Lemma 1

\begin{equation}\label{11}
\Delta^{odd}_{3,0}+\Delta^{odd}_{3,1}+\Delta^{odd}_{3,2}=\delta_n,
\end{equation}

where $|\delta_n|\leq 1$.

Thus by (\ref{3}) and (\ref{8})

\begin{equation}\label{12}
\Delta^{odd}_{3,2}=\delta_n-\Delta^{odd}_{3,0}-\Delta^{odd}_{3,1}=
\Delta^{odd}_3([0,n))-\Delta_3\left(\left[0,\frac n
2\right)\right)+\delta_n.
\end{equation}

With help of (\ref{12}) and the exact formulas for
$\Delta^{odd}_n(n),\;\;\Delta_3(n)$  \cite{4} we obtain in
particular that

\begin{equation}\label{13}
\Delta^{odd}_{3,2}([0,2^{2n-1}))=-3^{n-2},\;\;\;\Delta^{odd}_{3,2}([0,2^{2n}))=0.
\end{equation}

Nevertheless, it is sufficient for us to understand (\ref{3}) for
small $|\Delta^{odd}_{3,2}(n)|$ by the following way: if

\begin{equation}\label{14}
|\Delta^{odd}_{3,2}(n)|\leq \sqrt{n}\;\; then\;\;
|\Delta^{primes}_{3,2}(n)|=O\left(\frac{\sqrt{n}}{\ln{n}}\right).
\end{equation}

Now if $|\Delta^{odd}_{3,2}|> \sqrt{n}$ by (\ref{3}), (\ref{8}) and
(\ref{12}) we have

\begin{equation}\label{15}
\pi^0(n)-\pi^e(n)=\Delta^{primes}_{3,1}(n)+\Delta^{primes}_{3,2}(n)\sim
\frac{3\Delta^{odd}_3([0,n))}{\ln{n}}.
\end{equation}

Note that, according to \cite{4}

\begin{equation}\label{16}
\lim_{n\rightarrow\infty}\frac{\ln{\Delta^{odd}_3([0,n))}}{\ln{n}}=n^\alpha.
\end{equation}
\newpage

If $|\Delta^{odd}_{3,2}|\leq \sqrt{n}$ then by (\ref{3}), (\ref{8})
and (\ref{14}) we have

\begin{equation}\label{17}
\pi^0(n)-\pi^e(n)=\frac{3\Delta^{odd}_{3,1}(n)}{\ln{n}}(1+o(1))+O(\sqrt{n})\sim
\frac{3\Delta_3([0,\frac n 2))}{\ln{n}}.
\end{equation}

Now by (\ref{10}), (\ref{15})-(\ref{17}) we find as the final result
that

$$
\ln{(\pi^0(n)-\pi^e(n))}=\frac{\ln 3}{\ln 4}\ln n+o(\ln n)
$$
 and our Conjecture 2 follows.$\blacksquare$

      Note that from Conjecture 2 evidently follows the statement of
Conjecture 1 but only for sufficiently large $n\geq n_0$.
Unfortunately, until now we are not able to estimate $n_0$.

Note that by the way we obtain the limits

$$
\lim_{n\rightarrow\infty}\frac{\ln{(\pi^o_{3,1}(n)-\pi^e_{3,1}(n))}}{\ln{n}}=
\frac{\ln{3}}{\ln{4}};
$$

for $n_k=2^{2k-1},\;\;k\in\mathbb{N}$,

$$
\lim_{k\rightarrow\infty}\frac{\ln{(\pi^e_{3,2}(n_k)-\pi^o_{3,2}(n_k))}}{\ln{n_k}}=
\frac{\ln{3}}{\ln{4}}.
$$

\; \; \; \;\;\;\;\;

     The author is grateful to D.Berend, R.K.Guy, G.Martin and
T.D.Noe which show an interest in the considered conjectures.

\end{document}